# Polyhedral Cones of Magic Cubes and Squares

Maya Ahmed
Jesús De Loera
Raymond Hemmecke

*Dedicated to Eli Goodman and Ricky Pollack on the occasion of their 66th/67th birthdays.*

Abstract

Using computational algebraic geometry techniques and Hilbert bases of polyhedral cones we derive explicit formulas and generating functions for the number of magic squares and magic cubes.

Magic cubes and squares are very popular combinatorial objects (see [2, 19, 21] and their references). A *magic square* is a square matrix whose entries are nonnegative integers and whose row sums, column sums, and main diagonal sums add up to the same integer number $s$. We will call $s$ the *magic sum* of the square. In the literature there have been many variations on the definition of magic squares. For example, one popular variation of our definition adds the restriction of using the integers $1,\ldots,n^2$ as entries (such magic squares are commonly called *natural* or *pure* and a large part of the literature consists of procedures for constructing such examples, see [2, 19, 21]), but in this article the entries of the squares will be arbitrary nonnegative integers. We will consider other kinds of restrictions instead:

*Semi-magic squares* is the case when only the row and column sums are considered. This apparent simplification has in fact a very rich theory and several open questions remain (see [9, 14, 26] and references within. Semi-magic squares are called magic squares in these references). *Pandiagonal magic squares* are magic squares with the additional property that any broken-line diagonal sum adds up to the same integer (see Figure 1).

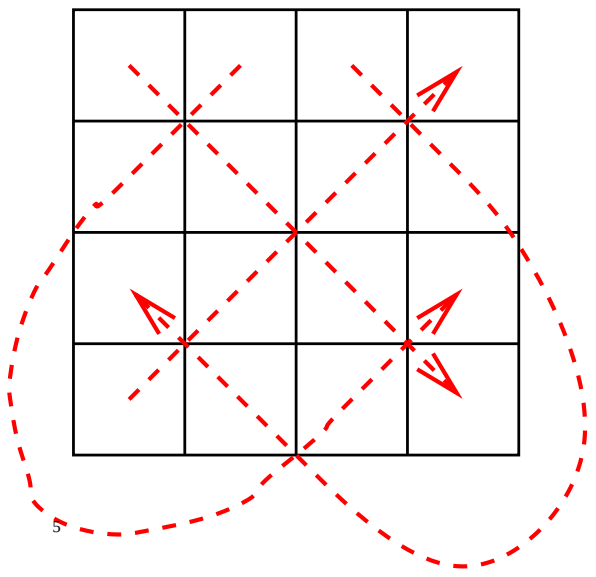

| | | | |
|---|---|---|---|
| 12 | 0 | 5 | 7 |
| 0 | 12 | 7 | 5 |
| 7 | 5 | 0 | 12 |
| 5 | 7 | 12 | 0 |

Figure 1: Four broken diagonals of a square and a pandiagonal magic square.

There are analogous definitions in higher dimensions. A *semi-magic hypercube* is a $d$-dimensional $n \times n \times \cdots \times n$ array of $n^d$ non-negative integers, which sum up to the same number $s$ for any line parallel to some axis. A *magic hypercube* is a semi-magic cube that has the additional property that the sums of all the main diagonals, the $2^{d-1}$ copies of the diagonal $x_{1,1,\ldots,1}, x_{2,2,\ldots,2}, \ldots, x_{n,n,\ldots,n}$ under the symmetries of the $d$-cube, are also equal to the magic sum. For example, in a $2 \times 2 \times 2$ cube there are 4 diagonals with sums $x_{1,1,1} + x_{2,2,2} = x_{2,1,1} + x_{1,2,2} = x_{1,1,2} + x_{2,2,1} = x_{1,2,1} + x_{2,1,2}$. We can see a magic

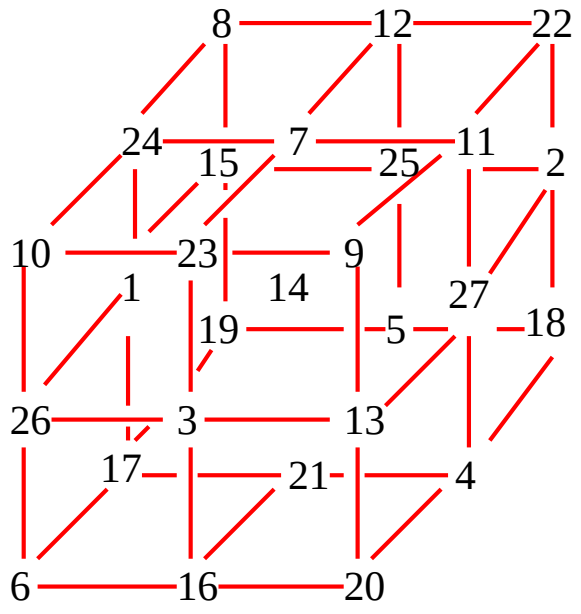

Figure 2: A magic cube.

Two fundamental problems about magic arrays are **(1)** enumerating such arrays and **(2)** generating particular elements. In this paper we address these two issues from a discrete-geometric perspective. The work of Ehrhart and Stanley [14, 15, 25, 26] when applied to the study of semi-magic squares showed that many enumerative and structural properties of magic arrays can actually be formulated in terms of polyhedral cones. The conditions of constant magic sum can be written in terms of a system $\{x|Ax = 0, x \geq 0\}$, where the vector $x$ has as many entries as there are cells in the array (labeled $x_{i_1,i_2,\ldots,i_d}$), and a matrix $A$ with entries $0, 1$ or $-1$ forces the different possible sums to be equal.

The purpose of this note is to study the convex polyhedral cones defined by magic squares, pandiagonal magic squares, semi-magic hypercubes, and magic hypercubes. In particular we study the Hilbert bases and extreme rays of these cones. We have used computational polyhedral geometry and commutative algebra techniques to derive explicit counting formulas for the four families of magic arrays we defined. Similar derivations had been done earlier for semi-magic squares [27, §4]. The interested reader can download the complete extreme ray information and Hilbert bases from `www.math.ucdavis.edu/~deloera/RESEARCH/magic.html`

Hilbert bases for these cones of magic arrays are special finite sets of nonnegative integer arrays that generate every other nonnegative integer array as a linear nonnegative integer combination of them. Most of our arguments will actually use minimal Hilbert bases which are smallest possible and unique [24]. Due to their size and complexity, our calculations of Hilbert bases and extreme rays were done with the help of a computer. We explain later on our algorithmic methods.

Having a Hilbert basis allows the generation of *any* magic array in the family, and makes trivial the construction of unlimited numbers of such objects or simply to list *all* magic arrays of fixed small size. Another benefit is that a Hilbert basis can be used to compute generating functions for the number of magic arrays from the computation of Hilbert series of the associated affine semigroup ring. We carry on these calculations using Gröbner bases methods. Finally minimal integer vectors along extreme rays of a cone are in fact also members of the Hilbert basis.

It is well-known from the work of Ehrhart [13] that for any rational pointed cone, if its lattice points receive a grading (e.g. by total sum of the entries, or in this case magic sum), then the function that counts the lattice points of fixed graded value is a *quasipolynomial*. A function $f : \mathbb{N} \to \mathbb{C}$ is quasipolynomial if there is an integer $N > 0$ and polynomials $f_0, \ldots, f_{N-1}$ such that $f(m) = f_i(m)$ if $m \equiv i$ $(mod\ N)$. The integer $N$, which is not unique, will be called a *quasi-period* of $f$. If it is the smallest quasi-period it will be called the *period of the quasipolynomial*. It is a natural question to investigate when the quasipolynomial is actually a polynomial, i.e when the period is one. We study this question for the four families of magic arrays. It is known that for the cone of semi-magic squares the quasi-polynomial is actually a polynomial (i.e. period is one). This follows from a well-known result of Ehrhart that assures that, for integral polytopes, the function that counts lattice points inside their integral dilations is a polynomial. We prove that the same argument does not work for other magic arrays. This will involve studying the extreme rays of the various cones. In general, determining exactly the period is a delicate issue as seen from Example 4.6.27 [27].

Consider the convex hull $P$ of real nonnegative arrays (of given size) all whose mandated sums equal 1. We call the polytope $P$ the *polytope of stochastic magic arrays*. For example, the stochastic semi-magic squares are the well-known bistochastic matrices ($n \times n$ matrices whose row and column sums are one) and $P$ is the famous Birkhoff-von Neumann polytope [9, 26]. It is easy to see that the polytope $P$ can be written as $P = \{x \in R^d : x \geq 0, \text{ and } Bx = 1\}$ where the matrix $B$ has $\{0, 1\}$ entries. $B$ has as many rows as axial sums (row, column, diagonals, etc), and the columns of $B$ correspond to the entries of the magic array.

In [6], Bona presented a proof that the counting function of semi-magic $3 \times 3 \times 3$ cubes is a quasi-polynomial of non-trivial period. In our first theorem we extend this by actually computing an explicit generating function and quasipolynomial formulas for the number of semi-magic $3 \times 3 \times 3$ cubes.

**Theorem 0.1.** *Denote by $SH_n^d(s)$ the number of semi-magic d-dimensional hypercubes with $n^d$ entries. We have the following results*

1. *From the Hilbert bases for the cones of $3\times3\times3$ semi-magic cubes we obtain the generating function.* $\sum_{s=0}^{\infty} SH_3^3(s)t^s = \frac{t^8+5t^7+67t^6+130t^5+242t^4+130t^3+67t^2+5t+1}{(1-t)^9(1+t)^2} (= 1 + 12t + 132t^2 + 847t^3 + 3921t^4 + 14286t^5 + 43687t^6 + 116757t^7 + \ldots)$.

   *In other words,*

$$SH_3^3(s) = \begin{cases} \frac{9}{2240}\,s^8 + \frac{27}{560}\,s^7 + \frac{87}{320}\,s^6 + \frac{297}{320}\,s^5 + \frac{1341}{640}\,s^4 + \frac{513}{160}\,s^3 + \frac{3653}{1120}\,s^2 + \frac{627}{280}\,s + 1 & \text{if } 2|s, \\ \frac{9}{2240}\,s^8 + \frac{27}{560}\,s^7 + \frac{87}{320}\,s^6 + \frac{297}{320}\,s^5 + \frac{1341}{640}\,s^4 + \frac{513}{160}\,s^3 + \frac{3653}{1120}\,s^2 + \frac{4071}{2240}\,s + \frac{47}{128} & \text{otherwise.} \end{cases}$$

2. *The number of vertices of the polytope of stochastic semi-magic $n \times n \times n$ cubes is bounded below by $(n!)^{2n}/n^{n^2}$. The polytopes of stochastic $3 \times 3 \times 3 \times 3$ semi-magic $3 \times 3 \times 3$ cubes and $3 \times 3 \times 3 \times 3$ hypercubes are not integral.*

We also computed an explicit generating function for the number of $3 \times 3 \times 3$ magic cubes.

**Theorem 0.2.** *Let $MC_n(s)$ denote the number of $n \times n \times n$ magic cubes. Then, $MC_n(s)$ is a quasipolynomial of degree $(n-1)^3 - 4$ for $n \geq 3, n \neq 4$. For $n = 4$ it has degree $(4-1)^3 - 3 = 24$. For $n = 3$, using the minimal Hilbert basis for the cones of $3 \times 3 \times 3$ magic cubes, we computed $\sum_{s=0}^{\infty} MC_3(s)t^s = \frac{t^{12}+14\,t^9+36\,t^6+14\,t^3+1}{(1-t^3)^5}$ $(= 1 + 19\,t^3 + 121\,t^6 + 439\,t^9 + 1171\,t^{12} + 2581\,t^{15} + 4999\,t^{18} + \ldots)$. Thus, in terms of a quasipolynomial formula we have:*

$$MC_3(s) = \begin{cases} \frac{11}{324}s^4 + \frac{11}{54}s^3 + \frac{25}{36}s^2 + \frac{7}{6}s + 1 & \text{if } 3|s, \\ 0 & \text{otherwise.} \end{cases}$$

*The polytope of stochastic $3 \times 3 \times 3 \times 3$ magic hypercubes is not integral.*

Our next contribution is to continue the enumerative analysis done in [4]. These authors wrote down formulas for the number of magic squares of orders 3 and 4. We have corrected a minor mistake in the $4 \times 4$ formula of [4, page 8] (the $3 \times 3$ case has been known since 1915 [20]), we find further values for order 5 magic squares and we give evidence supporting one of their conjectures [4, page 9].

**Theorem 0.3.** *If $M_n(s)$ denotes the number of $n \times n$ magic squares of magic sum $s$, then , from the minimal Hilbert bases for the cones of $4 \times 4$ and $5 \times 5$ magic squares, we obtain*

$\sum_{s=0}^{\infty} M_4(s)t^s = \frac{t^8+4t^7+18t^6+36t^5+50t^4+36t^3+18t^2+4t+1}{(1-t)^4(1-t^2)^4} (= 1+8t+48t^2+200t^3+675t^4+1904t^5+4736t^6+10608t^7+21925t^8+\dots)$*,*

*specifically we obtain that*

$$M_4(s) = \begin{cases} \frac{1}{480}s^7 + \frac{7}{240}s^6 + \frac{89}{480}s^5 + \frac{11}{16}s^4 + \frac{49}{30}s^3 + \frac{38}{15}s^2 + \frac{71}{30}s + 1 & \text{if } 2|s, \\ \frac{1}{480}s^7 + \frac{7}{240}s^6 + \frac{89}{480}s^5 + \frac{11}{16}s^4 + \frac{779}{480}s^3 + \frac{593}{240}s^2 + \frac{1051}{480}s + \frac{13}{16} & \text{otherwise.} \end{cases}$$

*We also know the values of $M_5(s)$ for $s \leq 6$. The polytope of stochastic magic squares is not integral for $n > 2$.*

Finally, we continue the work started in [1, 16] for the study of pandiagonal magic squares. Here we investigate their Hilbert bases, as an application we recomputed the formulas of Halleck (see [16, Chapters 8,10]. The integrality of the polytope of panstochastic magic squares was fully solved in [1].

**Theorem 0.4.** *Let $MP_n(s)$ denote the number of $n \times n$ pandiagonal magic squares of magic sum $s$, then from the Hilbert bases for the cones of $4 \times 4$ and $5 \times 5$ pandiagonal magic squares we obtain*

$$MP_4(s) = \begin{cases} \frac{1}{48}(s^2+4s+12)(s+2)^2 & \text{if } 2|s, \\ 0 & \text{otherwise.} \end{cases}$$

$$MP_5(s) = \frac{1}{8064}(s+4)(s+3)(s+2)(s+1)(s^2+5s+8)(s^2+5s+42).$$

Here is the plan for the paper: In Section 1 we review the notion of (minimal) Hilbert bases and how we computed them. We show how to use a Hilbert basis to compute a generating function that counts the number of nonnegative integer arrays of given magic sum. In that section we recall some basic facts about polyhedral cones, Ehrhart polynomials, and commutative semigroup rings (see [11, 28]). Finally, in Section 2, we discuss the specific details for the four theorems above, each appearing in a separate subsection. We close this introduction remarking that the algebraic-geometric techniques used here are not the only useful computational tools. In fact, there has been a surge of interest on such techniques with good practical results (see [3, 12, 29]).

# 1 Hilbert bases for counting and element generation

Let $A$ be an integer $d \times n$ matrix, we study pointed cones of the form $C = \{x | Ax = 0, x \geq 0\}$. A cone is *pointed*, if it does not contain any linear subspace besides the origin. It is well-known that pointed cones admit also a representation as the set of all possible nonnegative real linear combinations of finitely many vectors, the so called *extreme rays* of the cone (see page 232 of [24]). As an example we consider the cone of $3 \times 3$ magic matrices. This cone is defined by the system of equations

$$x_{11} + x_{12} + x_{13} = x_{21} + x_{22} + x_{23} = x_{31} + x_{32} + x_{33}$$

$$x_{11} + x_{12} + x_{13} = x_{11} + x_{21} + x_{31} = x_{12} + x_{22} + x_{32} = x_{13} + x_{23} + x_{33}$$

$$x_{11} + x_{12} + x_{13} = x_{11} + x_{22} + x_{33} = x_{31} + x_{22} + x_{13},$$

and the inequalities $x_{ij} \geq 0$. In our example for $3 \times 3$ magic squares the cone $C$ has dimension 3, it is a cone based on a quadrilateral, thus it has 4 rays (see Figure 3). It is easy to see that all other cones that we will treat for magic arrays are also solutions of a system $Ax = 0, x \geq 0$, where $A$ is a matrix with $0, 1, -1$ entries. For a given cone $C$ we are interested in $S_C = C \cap \mathbb{Z}^n$, the *semigroup of the cone* $C$.

An element $v$ of $S_C$ is called *irreducible* if a decomposition $v = v_1 + v_2$ for $v_1, v_2 \in S_C$ implies that $v_1 = 0$ or $v_2 = 0$. A *Hilbert basis* for $C$ is a finite set of vectors $HB(C)$ in $S_C$ such that every other element of $S_C$ is a positive integer combination of elements in $HB(C)$. A *minimal Hilbert basis* $HB(C)$ is inclusion minimal with respect to all other Hilbert bases of $C$. As a consequence all elements of the minimal Hilbert basis $HB(C)$ are irreducible and $HB(C)$ is unique.

| | | |
|---|---|---|
| 1 | 0 | 2 |
| 2 | 1 | 0 |
| 0 | 2 | 1 |
| 2 | 0 | 1 |
| 0 | 1 | 2 |
| 1 | 2 | 0 |
| 0 | 2 | 1 |
| 2 | 1 | 0 |
| 1 | 0 | 2 |
| 1 | 2 | 0 |
| 0 | 1 | 2 |
| 2 | 0 | 1 |
| 1 | 1 | 1 |
| 1 | 1 | 1 |
| 1 | 1 | 1 |

Figure 3: The Hilbert basis for the cone of $3 \times 3$ magic squares. The top four squares are the rays of the cone.

A natural question is then, *how can we compute the minimal Hilbert basis of a cone $C$?* Several research communities have developed algorithms for computing Hilbert bases having different applications in mind: integer programming and optimization [17], commutative algebra [7, 23, 28], and constraint programming [10, 22]. In our calculations of minimal Hilbert bases we used extensively the novel project-and-lift algorithm presented in [18] and implemented in `MLP` by R. Hemmecke. On the other hand we were able to corroborate independently most of our results using a different algorithm, the cone decomposition algorithm, implemented in `NORMALIZ` by Bruns and Koch [7]. Similar ideas were also discussed in [27]. Now we present brief descriptions of these two methods.

Hemmecke's algorithm for computing the Hilbert basis $H$ of a pointed rational cone $C$ expressed as $\{z : Az = 0, z \in \mathbb{R}^n_+\}$ proceeds as follows: Let $\pi_j : \mathbb{R}^n \to \mathbb{R}^j$ be the projection onto the first $j$ coordinates.

Let $K_j := \{\pi^n_j(v) : v \in \ker_{\mathbb{Z}}(A)\}$, $K^+_j := K_j \cap (\mathbb{R}^{j-1}_+ \times \mathbb{R}_+)$, and $K^-_j := K_j \cap (\mathbb{R}^{j-1}_+ \times \mathbb{R}_-)$, where $\ker_{\mathbb{Z}}(A) = \{z : Az = 0, z \in \mathbb{Z}^n\}$ denotes the integral kernel (or null space) of $A$. Observe that $K^+_j$ and $K^-_j$ are semi-groups under vector addition. Let $H^+_j$ and $H^-_j$ denote the unique inclusion minimal generating sets of the semi-groups $(K^+_j, +)$ and $(K^-_j, +)$. Clearly, $H = H^+_n$, since $K^+_n = C$.

The idea of the project-and-lift algorithm is to start with $H^+_1$, which is easy to compute, and to compute $H^+_{j+1} \cup H^-_{j+1}$ from $H^+_j$. This last step is done by a completion procedure (similar to s-pair reduction in Buchberger's algorithm [11]) and is based on the fact that for any vector $v \in \ker_{\mathbb{Z}}(A)$ with $\pi_{j+1}(v) \in H^+_{j+1} \cup H^-_{j+1}$, the vector $\pi_j(v)$ can be written as a non-negative integer linear combination of elements in $H^+_j$. Since many unnecessary vectors are already thrown away when $H^+_{j+1}$ is extracted from $H^+_{j+1} \cup H^-_{j+1}$, intermediate results are kept comparably small and larger problems can be solved.

The cone-decomposition algorithm, used in `NORMALIZ` triangulates the cone $C$ into finitely simplicial cones. A cone is simplicial if it is spanned by exactly $n$ linearly independent vectors $v_1, \ldots, v_n$. There are many possible triangulations, and any of these can be used. For each simplicial cone consider the parallepiped $\Pi = \{\lambda_1 v_1 + \cdots + \lambda_n v_n \in \mathbb{Z}^n | \lambda_i \in [0, 1), \}$. It is easy to see that the finite set of points $G_i = \Pi \cap \mathbb{Z}^n$ generates the semigroup. The computation of $G_i$ can be done via direct enumeration and the knowledge that $|G_i|$ is the same as the number of cosets of the quotient of $\mathbb{Z}^n$ by the Abelian group generated by the cone generators.

This way, each simplicial cone $\sigma_i$ in the triangulation of $C$ provides us with a set of generators $G_i$. From the union $G = \cup G_i = \{w_1, \ldots, w_m\}$, which obviously generates $C \cap \mathbb{Z}^n$, we need to find a subset $H \subset G$ whose elements are irreducible and still generate $C \cap \mathbb{Z}^n$. The subset $H$ is constructed recursively, starting from the empty set, in the $k$-th step we check if $w_k - h \in C$ for some $h \in H$. If yes, delete $w_k$ from the list and go to the next iteration; otherwise remove all those $h$ in $H$ which satisfy $h - w_k \in C$ and add $w_k$ to $H$ before passing to the next step. Clearly, since we have the inequality representation of the cone, it is easy to decide whether a vector belongs to the cone or not.

With any $d$-dimensional rational pointed polyhedral cone $C = \{Ax = 0, x \geq 0\}$ and a field $k$ we associate a *semigroup ring*, $R_C = k[y^a : a \in S_C]$, where there is one monomial $y_1^{a_1} y_2^{a_2} \ldots y_d^{a_d}$ in the ring for each element $a = (a_1, \ldots, a_d)$ of the semigroup $S_C$. By the definition of a Hilbert basis we know that every element of $S_C$ can be written as a finite linear combination $\sum \mu_i h_i$ where the $\mu_i$ are nonnegative integers. Thus $R_C$ is in fact a finitely generated $k$-algebra, with one generator per element of a Hilbert bases. Therefore $R_C$ can be written as the quotient $k[x_1, x_2, \ldots, x_N]/I_C$: Once we have the Hilbert basis $H = \{h_1, \ldots, h_N\}$ for the cone $C$, $I_C$ is simply the kernel of the polynomial map $\phi : k[x_1, x_2, \ldots, x_N] \longrightarrow k[y_1, y_2, \ldots y_d]$, where $\phi(x_i) = y^{h_i}$ and for $h_i = (a_1, a_2, \ldots, a_d)$ we set $y^{h_i} = y_1^{a_1} y_2^{a_2} \ldots y_d^{a_d}$. There are standard techniques for computing this kind of kernel (see [28] and references within).

It is important to observe that, for our cones of magic arrays, we can give a natural grading to $R_C$. A magic array can be thought of as a monomial on the ring and its degree will be its magic sum. For example, all the elements of the Hilbert bases of $3 \times 3$ magic squares are elements of degree 3. Once we have a graded $k$-algebra we can talk about its decomposition into the direct sum of its graded components $R_C = \bigoplus R_C(i)$, where each $R_C(i)$ collects all elements of degree $i$ and it is a $k$-vector space (where $R_C(0) = k$). The function $H(R_C, i) = dim_k(R_C(i))$ is the *Hilbert function* of $R_C$. Similarly one can construct the *Hilbert-Poincaré series* of $R_C$, $H_{R_C}(t) = \sum_{i=0}^{\infty} H(R_C, i) t^i$.

**Lemma 1.1.** *Let $C$ be a pointed rational cone, with Hilbert basis $H = \{h_1, \ldots, h_N\}$. Let the degree of a variable $x_i$ in the ring $k[x_1, \ldots, x_N]$ be the magic sum of the its corresponding Hilbert basis element $h_i$. Let $R_C$ be the (graded) semigroup ring obtained from the minimal Hilbert basis of a cone $C$ of magic arrays. Then the number of distinct magic arrays of magic constant $s$ equals the value of the Hilbert function $H(R_C, s)$.*

*Proof:* By the definition of a Hilbert basis we have that every magic array in the cone $C$ can be written as a linear integer combination of the elements of the Hilbert basis. The elements of $HB(C) = \{h_1, h_2, \ldots, h_N\}$ are not affinely independent therefore there are different combinations that produce the same magic array. We have some dependencies of the form $\sum a_i h_i = \sum a_j h_j$ where the sums run over some subsets of $\{1, \ldots, N\}$. We consider such identities as giving a single magic array. The dependencies are precisely the elements of the toric ideal $I_C$, that give $R_C = k[x_1, x_2, \ldots, x_N]/I_C$. Every such dependence is a linear combination of generators of any Gröbner basis of the ideal $I_C$. Thus, if we encode a magic array $X$ as a monomial in variables $x_1, \ldots, x_N$ whose exponents are the coefficients of the corresponding Hilbert basis elements that add to $X$, we are counting the equivalence classes modulo $I_C$. These are called *standard monomials.* Finally, it is known that the number of standard monomials of graded degree $i$ equals the dimension of $R_C(i)$ as a $k$-vector space [11, Chapter 9]. □

It is known that the Hilbert-Poincaré series of $R_C$ can be expressed as a rational function of the form $H_{R_C}(t) = \frac{p(t)}{\Pi_{i=1}^{r}(1-t^{\delta_i})}$. where $\delta_i$ can be read from the rays of the cone $C$; they correspond to the denominators of the vertices of the polytope of stochastic arrays whose dilations give the cone $C$ (see Theorem 4.6.25 [27] and Theorem 2.3 in [26]). To compute the Hilbert-Poincaré series we relied on the computer algebra package `CoCoA` [8], that has implementations for different algorithms of Hilbert series computations [5]. The basic idea comes from the theory of Gröbner bases (see [11, §9]). It is known that the initial ideal of $I_C$ with respect to any monomial order gives a monomial ideal $J$ and the Hilbert functions of $k[x_1, x_2, \ldots, x_N]/I_C$ and $k[x_1, x_2, \ldots, x_N]/J$ are equal. Computing the Hilbert function of the monomial ideal $J$ is a combinatorial problem which can be solved by an inclusion-exclusion type procedure [5] that eliminates variables at each iteration.

We illustrate the above algebraic techniques calculating a formula for the number of $3 \times 3$ magic squares, where $x_5$ corresponds to the matrix with all entries one, at the bottom of Figure 3, and the other 4 variables $x_1, x_2, x_3, x_4$ correspond to the magic squares on top of Figure 3, as they appear from left to right. The ideal $I_C$ given by the kernel of the map is generated by the two relations $x_1x_4 - x_5^2$, $x_2x_3 - x_1x_4$. The first relation means, for example, that the sum of magic square 1 with magic square 4 is the same as twice the magic square 5. The `CoCoA` commands that compute the Hilbert-Poincaré series is

```
L:=[3,3,3,3,3];
Use S::=Q[x[1..5]],Weights(L);
I:=Ideal(x[1]*x[4]-x[5]^2,x[1]*x[4]-x[2]*x[3]);
Poincare(S/I);
---  Non-simplified HilbertPoincare' Series  ---
(1 - 2x[1]^6 + x[1]^12) / ( (1-x[1]^3) (1-x[1]^3) (1-x[1]^3) (1-x[1]^3) (1-x[1]^3) )
```

Note that to carry out the computation it is necessary to specify a weight for the variables. In our case the weights are simply the magic sums of the array. It is known that from a rational representation like this one can directly recover a quasipolynomial (see [27, §4]).

$$M_3(s) = \begin{cases} \frac{2}{9}s^2 + \frac{2}{3}s + 1 & \text{if } 3|s, \\ 0 & \text{otherwise.} \end{cases}$$

We have seen already that magic arrays are nonnegative integer solutions of a system $Ax = 0, x \geq 0$, where $A$ is a matrix with $\{0, 1, -1\}$ entries. This system defines a pointed rational polyhedral cone $C$. One can set up the cone $C$ as the union of all (real-valued) dilations of the polytope of stochastic magic arrays $P = \{x \in R^d : x \geq 0, \text{ and } Bx = 1\}$. For a positive number $n$ we denote $E(P, n)$ the number of lattice points in the dilation $nP = \{nx | x \in P\} = \{x \in R^d : x \geq 0, \text{ and } Bx = n \cdot 1\}$. Note that when $n$ is integer, $E(P, n)$ for the $P$ polytope of stochastic magic arrays counts the number of integral magic arrays of magic sum $n$. If we let $n$ take real values, then the union of the different dilations of $P$ as $n$ changes is the pointed polyhedral cone $C$. This is easy to see since any magic square of magic sum $\lambda$ satisfies the equations $Ax = 0, x \geq 0$, thus all dilations are contained in the cone $C$. On the other hand, any solution $x$ to the system of inequalities that defines $C$ is a magic square of real valued magic sum $\lambda$. Dividing all entries of the array $x$ by $\lambda$ we obtained a magic array that satisfies the system $Bx = 1, x \geq 0$, thus the cone $C$ is contained in the union of all dilations of $P$. It can be verified that the rays of the cone $C$ are given by all scalar multiples of vertices of $P$. For our purposes the main result is a theorem of Ehrhart:

**Lemma 1.2 ([13]).** *For a rational $k$-polytope $P$ embedded in $\mathbb{R}^d$, in particular the polytope of stochastic magic arrays, the counting function $E(P, n)$ is a quasipolynomial in $n$ whose degree equals $k$ and whose period is less than or equal to the least common multiple of the denominators of the vertices of $P$.*

For example, for $3 \times 3$ magic squares the vertices of the polytope of stochastic magic squares are obtained by dividing the first 4 magic squares in Figure 3 by 3. In this case the periodicity of the function is exactly three. Although in all our computations the period of the quasipolynomial turned out to be equal to the least common multiple of the denominators of the vertices of $P$, this is not true in general (see Example 4.6.27 in [27]).

# 2 Families of Magic Arrays.

## 2.1 Semi-magic Hypercubes: Theorem 0.1

We consider first the $3 \times 3 \times 3$ semi-magic cube. Bona [6] had already observed that a Hilbert basis must contain only elements of magic constant one and two. Here we provide the 12 Hilbert basis elements of magic constant 1. There are 54 of magic constant 2, which we are not listing here, but can be downloaded from `www.math.ucdavis.edu/~deloera/RESEARCH/magic.html`

| | |
|---|---|
| (0, 0, 1, 1, 0, 0, 0, 1, 0, 1, 0, 0, 0, 1, 0, 0, 0, 1, 0, 1, 0, 0, 0, 1, 1, 0, 0) | (0, 0, 1, 0, 1, 0, 1, 0, 0, 0, 1, 0, 1, 0, 0, 0, 0, 1, 1, 0, 0, 0, 0, 1, 0, 1, 0) |
| (1, 0, 0, 0, 0, 1, 0, 1, 0, 0, 1, 0, 1, 0, 0, 0, 0, 1, 0, 0, 1, 0, 1, 0, 1, 0, 0) | (0, 0, 1, 1, 0, 0, 0, 1, 0, 0, 1, 0, 0, 0, 1, 1, 0, 0, 1, 0, 0, 0, 1, 0, 0, 0, 1) |
| (1, 0, 0, 0, 1, 0, 0, 0, 1, 0, 0, 1, 1, 0, 0, 0, 1, 0, 0, 1, 0, 0, 0, 1, 1, 0, 0) | (0, 0, 1, 0, 1, 0, 1, 0, 0, 1, 0, 0, 0, 0, 1, 0, 1, 0, 0, 1, 0, 1, 0, 0, 0, 0, 1) |
| (0, 1, 0, 1, 0, 0, 0, 0, 1, 0, 0, 1, 0, 1, 0, 1, 0, 0, 1, 0, 0, 0, 0, 1, 0, 1, 0) | (0, 1, 0, 1, 0, 0, 0, 0, 1, 1, 0, 0, 0, 0, 1, 0, 1, 0, 0, 0, 1, 0, 1, 0, 1, 0, 0) |
| (1, 0, 0, 0, 1, 0, 0, 0, 1, 0, 1, 0, 0, 0, 1, 1, 0, 0, 0, 0, 1, 1, 0, 0, 0, 1, 0) | (1, 0, 0, 0, 0, 1, 0, 1, 0, 0, 0, 1, 0, 1, 0, 1, 0, 0, 0, 1, 0, 1, 0, 0, 0, 0, 1) |
| (0, 1, 0, 0, 0, 1, 1, 0, 0, 1, 0, 0, 0, 1, 0, 0, 0, 1, 0, 0, 1, 1, 0, 0, 0, 1, 0) | (0, 1, 0, 0, 0, 1, 1, 0, 0, 0, 0, 1, 1, 0, 0, 0, 1, 0, 1, 0, 0, 0, 1, 0, 0, 0, 1) |

From the Hilbert basis and using `CoCoA` to compute the number of magic cubes we obtain the stated rational generating function. Now we claim the number of vertices of the polytope of stochastic semi-magic $n \times n \times n$ cubes is bounded below by $(n!)^{2n}/n^{n^2}$. This follows from a bijection between integral stochastic semi-magic cubes and $n \times n$ latin squares: Each 2-dimensional layer or slice of the integral stochastic cubes are permutation matrices (by Birkhoff-Von Neumann theorem), the different slices or layers cannot have overlapping entries else that would violate the fact that along a line the sum of the entries equals one. Thus make the permutation coming from the first slice be the first row of the latin square, the second slice permutation gives the second row of the latin square, etc. From well-known bounds for latin squares we obtain the lower bound.

The polytope of stochastic semi-magic $3 \times 3 \times 3$ cubes is actually not equal to the convex hull of integral semi-magic cubes. This follows for because the 54 elements of degree two in the Hilbert basis, when appropriately normalized, give rational stochastic matrices that are all vertices. In other words, the Birkhoff-von Neumann theorem [24, page 108] about stochastic semi-magic matrices is false for $3 \times 3 \times 3$ stochastic semi-magic cubes. Finally, the polytope of 4-dimensional semi-magic hypercubes has a non-integral vertex (each row is a 3-cube worth of values):

$$\begin{gathered}1/3 * (0,2,1,2,1,0,1,0,2,1,1,1,0,1,2,2,1,0,2,0,1,1,1,1,0,2,1\\ 2,0,1,0,1,2,1,2,0,1,2,0,1,1,1,1,0,2,0,1,2,2,1,0,1,1,1\\ 1,1,1,1,1,1,1,1,1,1,0,2,2,1,0,0,2,1,1,2,0,0,1,2,2,0,1)\end{gathered}$$

## 2.2 Magic Hypercubes: Theorem 0.2

The function that counts magic cubes is a quasipolynomial whose degree is the same as the dimension of the cone of magic cubes minus one. For small values (e.g $n = 3, 4$) we can directly compute this. We present an argument for its value for $n > 4$:

**Lemma 2.1.** *Let $B$ be the $(3n^2 + 4) \times n^3$ matrix with $0, 1$ entries determining axial and diagonal sums. In this way we see that $n \times n \times n$ magic cubes of magic sum $s$ are the integer solutions of $Bx = (s, s, \ldots, s)^T, x \geq 0$. For $n > 4$ the kernel of the matrix $B$ has dimension $(n-1)^3 - 4$.*

**Proof.** It is known that for semi-magic cubes the dimension is $(n-1)^3$ [4], which means that the rank of the submatrix $B'$ of $B$ without the 4 rows that state diagonal sums is $n^3 - (n-1)^3$. It remains to be shown that the addition of the 4 sum constraints on the main diagonals to the defining equations of the $n \times n \times n$ semi-magic cube increases the rank of the defining matrix $B$ by exactly 4.

Let us denote the $n^3$ entries of the cube by $x_{1,1,1}, \ldots, x_{n,n,n}$ and consider the $(n-1)\times(n-1)\times(n-1)$ sub-cube with entries $x_{1,1,1}, \ldots, x_{n-1,n-1,n-1}$. For a semi-magic cube we have complete freedom to choose these $(n-1)^3$ entries. The remaining entries of the $n \times n \times n$ magic cube become known via the semi-magic cube equations, and all entries together form a semi-magic cube. For example:

$$x_{n,1,1} = -\sum_{i=1}^{n-1} x_{i,1,1},\; x_{1,n,n} = \sum_{i=1}^{n-1}\sum_{j=1}^{n-1} x_{i,j,1},\; x_{n,n,n} = -\sum_{i=1}^{n-1}\sum_{j=1}^{n-1}\sum_{k=1}^{n-1} x_{i,j,k}.$$

However, for the magic cube, 4 more conditions have to be satisfied along the main diagonals. Employing the above semi-magic cube equations, we can rewrite these 4 equations for the main diagonals such that they involve only the variables $x_{1,1,1}, \ldots, x_{n-1,n-1,n-1}$. Thus, as we will see, the complete freedom of choosing values for the variables $x_{1,1,1}, \ldots, x_{n-1,n-1,n-1}$ is restricted by 4 independent equations. Therefore the dimension of the kernel of $B$ is reduced by 4.

Let us consider the 3 equations in $x_{1,1,1}, \ldots, x_{n-1,n-1,n-1}$ corresponding to the main diagonals $x_{1,1,n}, \ldots, x_{n,n,1}$, $x_{1,n,1}, \ldots, x_{n,1,n}$, and $x_{n,1,1}, \ldots, x_{1,n,n}$. They are linearly independent, since the variables $x_{n-1,n-1,1}$, $x_{n-1,1,n-1}$, and $x_{1,n-1,n-1}$ appear in exactly one of these equations. The equation corresponding to the diagonal $x_{1,1,1}, \ldots, x_{n,n,n}$ is linearly independent from the other 3, because, when rewritten in terms of only variables of the form $x_{i,j,k}$ with $1 \leq i,j,k < n$, it contains the variable $x_{2,2,3}$, which for $n > 4$ does not lie on a main diagonal and is therefore not involved in one of the other 3 equations. This completes the proof. □

We consider now the $3 \times 3 \times 3$ magic cubes. There are 19 elements in the Hilbert basis and all of them have magic sum value of 3. This already indicates that there is a quasipolynomial counting formula since there are no elements of magic sum not divisible by 3.

$(2,1,0,1,0,2,0,2,1,0,2,1,2,1,0,1,0,2,1,0,2,0,2,1,2,1,0)$
$(1,2,0,1,0,2,1,1,1,1,0,2,2,1,0,0,2,1,1,1,1,0,2,1,2,0,1)$
$(2,0,1,0,2,1,1,1,1,0,2,1,2,1,0,1,0,2,1,1,1,1,0,2,1,2,0)$
$(1,1,1,2,0,1,0,2,1,1,2,0,0,1,2,2,0,1,1,0,2,1,2,0,1,1,1)$
$(1,2,0,2,0,1,0,1,2,2,0,1,0,1,2,1,2,0,0,1,2,1,2,0,2,0,1)$
$(1,1,1,1,1,1,1,1,1,1,1,1,1,1,1,1,1,1,1,1,1,1,1,1,1,1,1)$
$(1,0,2,0,2,1,2,1,0,0,2,1,2,1,0,1,0,2,2,1,0,1,0,2,0,2,1)$
$(1,1,1,0,2,1,2,0,1,1,0,2,2,1,0,0,2,1,1,2,0,1,0,2,1,1,1)$
$(1,0,2,1,2,0,1,1,1,1,2,0,0,1,2,2,0,1,1,1,1,2,0,1,0,2,1)$
$(0,1,2,1,2,0,2,0,1,2,0,1,0,1,2,1,2,0,1,2,0,2,0,1,0,1,2)$
$(1,1,1,1,0,2,1,2,0,0,2,1,2,1,0,1,0,2,2,0,1,0,2,1,1,1,1)$
$(2,1,0,1,1,1,0,1,2,1,1,1,1,1,1,1,1,1,0,1,2,1,1,1,2,1,0)$
$(2,1,0,0,2,1,1,0,2,1,0,2,2,1,0,0,2,1,0,2,1,1,0,2,2,1,0)$
$(0,1,2,2,0,1,1,2,0,1,2,0,0,1,2,2,0,1,2,0,1,1,2,0,0,1,2)$
$(0,2,1,1,0,2,2,1,0,1,0,2,2,1,0,0,2,1,2,1,0,0,2,1,1,0,2)$
$(2,0,1,1,2,0,0,1,2,1,2,0,0,1,2,2,0,1,0,1,2,2,0,1,1,2,0)$
$(0,2,1,2,0,1,1,1,1,2,0,1,0,1,2,1,2,0,1,1,1,1,2,0,1,0,2)$
$(0,1,2,1,1,1,2,1,0,1,1,1,1,1,1,1,1,1,2,1,0,1,1,1,0,1,2)$
$(1,1,1,1,2,0,1,0,2,2,0,1,0,1,2,1,2,0,0,2,1,2,0,1,1,1,1)$

From this information, and using `CoCoA`, we can derive the desired formula for the count that appears in Theorem 0.2. Finally we include below an extreme ray for the cone of magic $3 \times 3 \times 3 \times 3$ hypercubes. Dividing its entries by 15 we get a rational vertex of the polytope of stochastic magic $3 \times 3 \times 3 \times 3$ hypercubes.

| | | | | | | | | | | | | | | | | | | | | | | | | | | |
|---|---|---|---|---|---|---|---|---|---|---|---|---|---|---|---|---|---|---|---|---|---|---|---|---|---|---|
| 8 | 7 | 0 | 0 | 8 | 7 | 7 | 0 | 8 | 4 | 4 | 7 | 5 | 2 | 8 | 6 | 9 | 0 | 3 | 4 | 8 | 10 | 5 | 0 | 2 | 6 | 7 |
| 4 | 4 | 7 | 5 | 2 | 8 | 6 | 9 | 0 | 1 | 10 | 4 | 8 | 5 | 2 | 6 | 0 | 9 | 10 | 1 | 4 | 2 | 8 | 5 | 3 | 6 | 6 |
| 3 | 4 | 8 | 10 | 5 | 0 | 2 | 6 | 7 | 10 | 1 | 4 | 2 | 8 | 5 | 3 | 6 | 6 | 2 | 10 | 3 | 3 | 2 | 10 | 10 | 3 | 2 |

## 2.3 Magic Squares: Theorem 0.3

$4 \times 4$ **magic squares:** Our calculations using `MLP` show that there are 20 elements in the Hilbert basis for the cone $C_{M4\times 4}$ of $4 \times 4$ magic squares. The 8 elements of magic sum one (not 7 as reported in [4]) and the 12 elements of magic sum 2 are listed below. To save space we present the squares as vectors $(x_{11}, \ldots, x_{14}, x_{21}, \ldots, x_{24}, x_{31}, \ldots, x_{34}, x_{14}, \ldots, x_{44})$.

$$\begin{array}{lll}
(0,0,1,0,0,1,0,0,0,0,0,1,1,0,0,0) & (0,0,1,0,1,0,0,0,0,1,0,0,0,0,0,1) & (1,0,0,0,0,0,1,0,0,0,0,1,0,1,0,0) \\
(0,0,0,1,1,0,0,0,0,0,1,0,0,1,0,0) & (0,1,0,0,0,0,0,1,0,0,1,0,1,0,0,0) & (1,0,0,0,0,0,0,1,0,1,0,0,0,0,1,0) \\
(0,0,0,1,0,1,0,0,1,0,0,0,0,0,1,0) & (0,1,0,0,0,0,1,0,1,0,0,0,0,0,0,1) &
\end{array}$$

These 8 permutation matrices are exactly all the magic squares of magic sum 1. The rest of the minimal Hilbert basis consists of magic sum 2 magic squares:

$$\begin{array}{lll}
(1,0,1,0,0,0,0,2,0,1,1,0,1,1,0,0) & (0,0,2,0,0,1,0,1,1,1,0,0,1,0,0,1) & (0,0,1,1,0,1,0,1,1,0,1,0,1,1,0,0) \\
(1,1,0,0,0,1,1,0,0,0,0,2,1,0,1,0) & (1,0,0,1,1,1,0,0,0,1,0,1,0,0,2,0) & (0,1,0,1,1,1,0,0,0,0,1,1,1,0,1,0) \\
(1,0,0,1,0,0,1,1,1,0,1,0,0,2,0,0) & (0,2,0,0,1,0,1,0,0,0,1,1,1,0,0,1) & (1,1,0,0,1,0,1,0,0,1,0,1,0,0,1,1) \\
(0,1,0,1,2,0,0,0,0,1,1,0,0,0,1,1) & (1,0,1,0,0,0,1,1,1,1,0,0,0,1,0,1) & (0,0,1,1,0,1,1,0,2,0,0,0,0,1,0,1)
\end{array}$$

Using `CoCoA`'s Hilbert series computation we obtain the generating function stated in Theorem 0.3.

$5 \times 5$ **magic squares** The $5 \times 5$ magic squares are the first challenging case. We were unable so far to recover the Hilbert series for this case. By using the fact that the Hilbert basis is a generating set we can easily compute several values of the Hilbert function, i.e. the numbers of magic squares for small values of the magic sum. Using the generators we consider all possible sums of them with small coefficients, making sure that repeated squares are only counted once. The values below allow us to prove that there is no polynomial formula that fits those values via interpolation. We use the Ehrhart-Macdonald reciprocity laws [26] that give us other 6 values of the function, which together with known roots allow for interpolation. There is no solution for the resulting linear system.

| magic sum | total number magic squares |
|---|---|
| 1 | 20 |
| 2 | 449 |
| 3 | 6792 |
| 4 | 67,063 |
| 5 | 484,419 |
| 6 | 2,750,715 |

The following table lists the number of elements in the Hilbert basis. All the elements for all the Hilbert bases we have computed can be obtained at `www.math.ucdavis.edu/~deloera/RESEARCH/magic.html`

| magic sum | number of HB elements |
|---|---|
| 1 | 20 |
| 2 | 240 |
| 3 | 1392 |
| 4 | 1584 |
| 5 | 1192 |
| 6 | 160 |
| 7 | 224 |
| 9 | 16 |
| | 4828 |

Finally we prove the rest of Theorem 0.3. We construct integral extreme ray vectors that, when its entries are divided by 1/2, give a fractional vertex of the polytope of stochastic magic squares: Let $n \geq 6$ and let $P_{n-2}$ be an $(n-2) \times (n-2)$ permutation matrix that does not contain a non-zero entry on its two main diagonals. Let $R_n$ be the $n \times n$ matrix that is constructed as follows:

- $R_{n,i,j} = 2 * P_{n-2,i-1,j-1}$ for $i = 2, \ldots, n-1$, $j = 2, \ldots, n-1$,
- $R_{n,1,j} = R_{n,n,j} = 0$ for $j = 2, \ldots, n-1$,
- $R_{n,i,1} = R_{n,i,n} = 0$ for $i = 2, \ldots, n-1$,
- $R_{n,1,1} = R_{n,n,1} = R_{n,1,n} = R_{n,n,n} = 1$.

Since $n-2 \geq 4$, there exists a permutation matrix $P_{n-2}$ with no non-zero entries on its main diagonals. Thus, $R_n$ is well-defined.

**Lemma 2.2.** *By construction, $R_n$ is a magic square of size $n$ with magic constant 2, in addition, for $n \geq 6$, $R_n$ is an extremal ray of the cone of $n \times n$ magic squares.*

*Proof:* Suppose that $R_n$ is not an extremal ray of the magic square cone. Therefore, there exists a non-zero magic square $\bar{R}_n$ with magic constant $s > 0$ whose support is strictly contained in the support of $R_n$. Since every row and column must have at least one non-zero entry, $\bar{R}_n$ must have a zero in one of the corners, that is without loss of generality, $\bar{R}_{n,1,1} = 0$. Since $s = \sum_{i=1}^{n} \bar{R}_{n,i,1} = \sum_{j=1}^{n} \bar{R}_{n,n,j}$, we obtain $s = \bar{R}_{n,n,1} = \bar{R}_{n,n,1} + \bar{R}_{n,n,n}$. Thus, $\bar{R}_{n,n,n} = 0$. But this contradicts $0 < s = \sum_{i=1}^{n} \bar{R}_{n,i,i} = 0$. Therefore, $\bar{R}_n$ does not exist, implying that $R_n$ is an extremal ray. $\square$

## 2.4 Pandiagonal Magic Squares: Theorem 0.4

Let us denote by $MP_n(s)$ the number of $n \times n$ pandiagonal magic squares with magic sum $s$. As in the case of magic squares the function $MP_n(s)$ is a quasipolynomial in $s$ of degree equal to the dimension of the cone plus one. Halleck [16] computed the dimension of the cone to be $(n-2)^2$ for odd $n$ and $(n-2)^2 + 1$ for even $n$ (degree of the quasipolynomial $MP_n(s)$ is one less than these). For the $4 \times 4$ pandiagonal magic squares a fast calculation corroborates that there are 8, magic-sum-2, generators. In his investigations, Halleck [16] identified a much larger generating set.

| | | |
|---|---|---|
| (1, 0, 0, 1, 0, 1, 1, 0, 1, 0, 0, 1, 0, 1, 1, 0) | (1, 0, 1, 0, 0, 1, 0, 1, 0, 1, 0, 1, 1, 0, 1, 0) | (0, 0, 1, 1, 1, 1, 0, 0, 0, 0, 1, 1, 1, 1, 0, 0) |
| (1, 0, 1, 0, 1, 0, 1, 0, 0, 1, 0, 1, 0, 1, 0, 1) | (0, 1, 0, 1, 0, 1, 0, 1, 1, 0, 1, 0, 1, 0, 1, 0) | (0, 1, 1, 0, 1, 0, 0, 1, 0, 1, 1, 0, 1, 0, 0, 1) |
| (1, 1, 0, 0, 0, 0, 1, 1, 1, 1, 0, 0, 0, 0, 1, 1) | (0, 1, 0, 1, 1, 0, 1, 0, 1, 0, 1, 0, 0, 1, 0, 1) | |

From the Hilbert basis we can calculate the formula stated in Theorem 0.4 using `CoCoA`. Finally we verify that the $5 \times 5$ pandiagonal magic squares have indeed a polynomial counting formula. This case requires in fact no calculations thanks to earlier work by [1] who proved that for $n = 5$ the only pandiagonal rays are precisely the pandiagonal permutation matrices. It is easy to see that only 10 of the 120 permutation matrices of order 5 are pandiagonal:

| | |
|---|---|
| (0, 0, 0, 1, 0, 1, 0, 0, 0, 0, 0, 0, 1, 0, 0, 0, 0, 0, 1, 0, 1, 0, 0, 0) | (1, 0, 0, 0, 0, 0, 0, 0, 1, 0, 0, 1, 0, 0, 0, 0, 0, 0, 0, 1, 0, 0, 1, 0, 0) |
| (0, 0, 1, 0, 0, 1, 0, 0, 0, 0, 0, 0, 0, 1, 0, 0, 1, 0, 0, 0, 0, 0, 0, 0, 1) | (0, 0, 0, 0, 1, 0, 1, 0, 0, 0, 0, 0, 1, 0, 1, 0, 0, 0, 0, 0, 0, 1, 0, 0) |
| (0, 0, 1, 0, 0, 0, 0, 0, 0, 1, 0, 1, 0, 0, 0, 0, 0, 0, 1, 0, 1, 0, 0, 0, 0) | (0, 0, 0, 0, 1, 0, 0, 1, 0, 0, 1, 0, 0, 0, 0, 0, 0, 0, 1, 0, 0, 1, 0, 0, 0) |
| (0, 0, 0, 1, 0, 0, 1, 0, 0, 0, 0, 0, 0, 0, 1, 0, 0, 1, 0, 0, 1, 0, 0, 0, 0) | (1, 0, 0, 0, 0, 0, 0, 1, 0, 0, 0, 0, 0, 0, 1, 0, 1, 0, 0, 0, 0, 0, 0, 1, 0) |
| (0, 1, 0, 0, 0, 0, 0, 0, 0, 1, 0, 0, 1, 0, 0, 1, 0, 0, 0, 0, 0, 0, 0, 1, 0) | (0, 1, 0, 0, 0, 0, 0, 0, 1, 0, 1, 0, 0, 0, 0, 0, 0, 1, 0, 0, 0, 0, 0, 0, 1) |

Once more a simple `CoCoA` calculation shows that the counting function equals indeed a polynomial, $\frac{1}{8064}(s+4)(s+3)(s+2)(s+1)(s^2+5s+8)(s^2+5s+42)$, as claimed in the Theorem.

## About Authors

Maya Ahmed, Jesús De Loera, and Raymond Hemmecke are at the Department of Mathematics, University of California, Davis, CA, USA, ahmed@math.ucdavis.edu, deloera@math.ucdavis.edu, raymond@hemmecke.de. This research was supported by NSF Grant DMS-0073815.